\documentclass[twoside]{ima}

\usepackage{amsmath}

\usepackage{amssymb}

\usepackage{graphicx,times}

\def\R{\mathbb{R}}
\def\RP{\mathbb{RP}}

\def\hom{\mathrm{Hom}}
\def\dim{\mathrm{dim\,}}
\def\ker{\mathrm{Ker\,}}
\def\im{\mathrm{Im\,}}
\def\span{\mathrm{span\,}}
\def\A{\mathbf{A}}
\def\B{\mathbf{B}}

\def\M{\mathbf{M}}
\def\g{\mathfrak{g}}
\def\a{\mathfrak{a}}
\def\b{\mathfrak{b}}

\def\m{\mathfrak{m}}

\def\sl{\mathfrak{sl}}
\def\so{\mathfrak{so}}

\def\SO{\mathrm{SO}}

\def\I{\mathcal{I}}

\begin{document}

\markboth{E. MUSSO, L. NICOLODI}{DIFFERENTIAL SYSTEMS AND TABLEAUX
OVER LIE ALGEBRAS}

\title{DIFFERENTIAL SYSTEMS ASSOCIATED WITH \\ TABLEAUX OVER LIE ALGEBRAS}

\author{Emilio Musso\thanks{
Dipartimento di Matematica Pura ed Applicata, Universit\`a degli
Studi dell'Aquila, Via Vetoio, I-67010 Coppito (L'Aquila), Italy
(musso@univaq.it). Partially supported by MIUR project {\it Metriche
riemanniane e variet\`a differenziali}, and by the GNSAGA of INDAM.}
\and Lorenzo Nicolodi\thanks{ Di\-par\-ti\-men\-to di
Ma\-te\-ma\-ti\-ca, Uni\-ver\-si\-t\`a degli Studi di Parma, Viale
G. P. Usberti 53/A - Campus universitario, I-43100 Parma, Italy}
(lorenzo.nicolodi@unipr.it). Partially supported by MIUR project
{\it Propriet\`a geometriche delle variet\`a reali e complesse}, and
by the GNSAGA of INDAM.}

\maketitle

\begin{abstract}
We give an account of the construction of exterior differential
systems based on the notion of tableaux over Lie algebras as
developed in \cite{MNtableaux}. The definition of a tableau over a
Lie algebra is revisited and extended in the light of the formalism
of the Spencer cohomology; the question of involutiveness for the
associated systems and their prolongations is addressed; examples
are discussed.
\end{abstract}


{\AMSMOS
  58A17, 58A15
 \endAMSMOS}

\section{Introduction}\label{s:intro}
The search for a common structure to various exterior differential
systems (EDSs) of geometric and analytic origin led to the algebraic
notion of a \textit{tableaux over a Lie algebra} \cite{MNtableaux}.
This notion builds on that of involutive tableau in the theory of
EDSs and can be seen as a non-commutative generalization of it.
Interestingly enough, from a tableau over a Lie algebra we can
canonically construct a linear Pfaffian differential system (PDS)
which is in involution and whose Cartan characters coincide with the
characters of the tableau.

\vskip0.1cm
Particular cases of this scheme lead to differential systems describing
well-known integrable systems such as the Grassmannian systems of Terng
\cite{Ter97, BDPT02}, the curved flat system
of Ferus and Pedit \cite{FPcurvedflats}, and many integrable surfaces
arising in projective differential geometry \cite{AG, Fe2, Fi}.
The tableaux corresponding to Grassmannian and curved
flat systems, the so-called Cartan tableaux, are obtained from the
Cartan decompositions of semisimple Lie algebras.
The tableaux corresponding to the various classes of integrable surfaces
in projective 3-space are given by sub-tableaux of a special
tableau over $\sl (4,\R)$. This is constructed by the method of moving
frames and amounts to the construction of a canonical adapted frame along
a generic surface in projective space (e.g., the Wilczynski--Cartan frame
\cite{Car1920defo, Fi, Bol, Fe2, Fe3}).
The involutiveness of these examples follows from the involutiveness
of the corresponding tableaux. The result that the Grassmannian system
and the curved flat system are in involution was first proved by Bryant
in a cycle of seminars at MSRI \cite{Br} and was then taken up and further
elaborated by Terng and Wang \cite{TW03}.

\vskip0.1cm

Other examples of linear Pfaffian system in involution which fit
into the above scheme, and in fact motivated our work, include: the
differential systems of isothermic surfaces in M\"obius and Laguerre
geometry \cite{Mtrieste, MNbudapest, MNbollettino, MNvarsavia}; the
differential systems of M\"obius-minimal (M-minimal or Willmore)
surfaces and of Laguerre-minimal (L-minimal) surfaces
\cite{MNtransactions, MNtrieste, Mtrieste}; the differential systems
associated to the deformation problem in projective geometry and in
Lie sphere geometry \cite{Car1920defo, Car-strasbourg-defo,
Blaschke1929, Jensen81, MNtohoku, Fe1}. The tableaux associated with
all these examples are constructed by the method of moving frames on
the Lie algebras of the corresponding symmetry groups (cf. Section
\ref{ss:conf-ex}). If, on the one hand, the presented approach may
be seen as a possibility to discuss various involutive systems from
a unified viewpoint, on the other hand, it can be viewed as a
possibility to find new classes of involutive systems. In this
respect, we mention the class of projective surfaces introduced in
\cite{MNtableaux}, which generalizes asymptotically-isothermic
surfaces and surfaces with constant curvature of Fubini's quadratic
form. An analogous class of surfaces in the context of conformal
geometry is discussed in Section \ref{ss:conf-ex}. For an
application of the above construction to the study of the Cauchy
problem for the associated systems, we refer the reader to
\cite{MNjmp, MNphysD}.

\vskip0.1cm

In this article we revisit the definition of tableau over a Lie
algebra using the formalism of the Spencer cohomology and extend it
to include 2-acyclic tableaux. This allows also non-involutive
systems into the scheme, reducing the question of involutiveness of
their prolonged systems to that of the prolongations of the
associated tableaux (cf. Section \ref{tableau-system}).

\vskip0.1cm

Section \ref{s:tableaux} contains the basic material about tableaux.
Section \ref{s:lie-tableaux} introduces the notion of tableaux over
Lie algebras. Section \ref{tableau-system} presents the construction
of EDSs from tableaux over Lie algebras and discusses some
properties of such systems. Section \ref{s:examples} discusses some
examples as an illustration of the theory developed in the previous
sections. Further developments are indicated in Section
\ref{s:developments}. The appendix collects some facts about the
Spencer complex of a tableau and the torsion of a PDS.

\section{Tableaux}\label{s:tableaux}

In this section we provide a summary of the results in the algebraic
theory of tableaux. As basic reference, we use the book by Bryant,
\textit{et al.} \cite{BCGGG}. See also \cite{IveyLandsberg2003}.

\vskip0.2cm
A \textbf{tableau} is a linear subspace $\A \subset\hom (\a,\b)$,
where $\a$, $\b$ are (real or complex) finite dimensional vector spaces.

An $h$-dimensional subspace $\a_h \subset \a$ is called \textit{generic}
w.r.t. $\A$ if the dimension of
\[
  \ker (\A,\a_h) : =\{Q\in \A \,|\, Q_{|_{\a_h}} =0\}
    \]
is a minimum.
The set of $h$-dimensional generic subspaces is a Zariski open
of the Grassmannian of $h$-dimensional subspaces of $\a$.

A flag $(0) \subset \a_1\subset \dots \subset \a_n=\a$ of $\a$
is said \textit{generic} if $\a_h$ is generic, for all
$h=1,\dots,n$.

The \textbf{characters} of $\A$ are the non-negative integers
$s_j(\mathbf{A})$, $j=1,\dots,n$, defined inductively by
\[
   s_1(\A) + \cdots +s_j(\A)= \mathrm{codim\,}
    \ker(\A,\a_j),
      \]
for any generic flag $(0)\subset \a_1\subset\cdots\subset \a_n=\a$.

From the definition, it is clear that
\[
 \dim \b \geq s_1\geq
     s_2\geq \dots \geq s_n,\quad \dim\A = s_1+\cdots +s_n.
     \]
If $s_{\nu}\neq 0$, but $s_{\nu +1}=0$, we say that $\A$ has
\textit{principal character} $s_{\nu}$ and call $\nu$ the
\textit{Cartan integer} of $\A$.

The \textbf{first prolongation} ${\A}^{(1)}$ of $\A$ is the kernel
of the linear map (Spencer coboundary operator, cf. Appendix
\ref{ss:spencer})
\[
 \delta^{1,1} : \hom(\a,\A)\cong \A\otimes \a^\ast \to \b\otimes \Lambda^2(\a^\ast)
  \]
\[
 \delta^{1,1}(F) (A_1,A_2) := \frac{1}{2}\left(F(A_1)(A_2) - F(A_2)(A_1)\right),
  \]
for $F\in \hom(\a,\A)$, and $A_1,A_2 \in \a$.

The \textbf{$h$-th prolongation} of ${\A}$ is defined inductively by
setting
\[
  {\A}^{(h)} = {{\A}^{(h-1)}}^{(1)},
   \]
for $h\geq 1$ (by convention ${\A}^{(0)}= \A$ and
${\A}^{(-1)}= \b$).
${\A}^{(h)}$ identifies with
\[
  \mathbf{A}^{(h)}=\left(\mathbf{A}\otimes S^h(\a^*)\right)\cap
     \left(\b\otimes S^{h+1}(\a^\ast)\right).
       \]
An element $Q_{(h)}\in \b\otimes S^{h+1}(\a^\ast)$ belongs to ${\A}^{(h)}$
if and only if $i(X)Q_{(h)} \in {\A}^{(h-1)}$, for all $X \in \a$.

\begin{theorem}
 $\dim {\A}^{(1)} \leq s_1 +2s_2 + \cdots + ns_n.$
   \end{theorem}

\vskip0.1cm
\begin{definition}
$\A$ is said \textbf{involutive} (or \textbf{in involution}) if
equality holds in the previous inequality.
\end{definition}

\begin{theorem}
For any tableau $\A$ there exists an integer $h_0$ such that
${\A}^{(h)}$ is involutive, for all $h\geq h_0$.
\end{theorem}

\begin{theorem}
If $\A$ is involutive, then ${\A}^{(1)}$ is involutive and
$$
  s^{(1)}_j: = s_j({\A}^{(1)}) = s_n(\A)+\cdots + s_j(\A),
     $$
$j=1,\dots,n$.
\end{theorem}

Thus \textit{every prolongation of an involutive tableau is involutive}.
Moreover, \textit{the principal character and
the Cartan integer are invariant under prolongation of an involutive tableau}.

It is well-known that $\A$ is involutive if and only if
\[
 H^{q,p}(\A)=(0),\quad \text{for all} \quad q\geq 1, p\geq 0
  \]
(Guillemin-Sternberg, Serre \cite{GuiSte64}). See Appendix
\ref{ss:spencer} for the definition of the Spencer groups
$H^{q,p}(\A)$.

\vskip0.2cm

A weaker notion is the following.

\begin{definition}
A tableau $\A$ is said
2-\textbf{acyclic} if
\[
  H^{q,2}(\A)=(0), \quad \text{for all} \quad q\ge 1.
   \]
\end{definition}

This notion plays an essential role in the prolongation procedure of
a non-involutive linear Pfaffian system (cf. Kuranishi, Goldschmidt
\cite{Goldschmidt1967a, Goldschmidt1967b}).

\vskip0.1cm

As shown in the following examples, tableaux and their prolongations
arise naturally in the context of PDE systems and of exterior
differential systems.

\vskip0.1cm

\begin{example}
{\rm
Let $V$ and $W$ be vector spaces with coordinates $x^1,\dots,x^n$
and $y^1,\dots,y^s$ dual to bases $v_1,\dots,v_n$ for $V$ and
$w_1,\dots,w_s$ for $W$. Consider the first-order constant
coefficient system of PDEs for maps $f : V \to W$ given in
coordinates by
\begin{equation}\label{lin-pde}
 B^{\lambda i}_{a}\,
 \frac{\partial y^a}{\partial{x^i}}(x) = 0 \quad (\lambda = 1,\dots,r).
  \end{equation}
The linear solutions $y^a(x) = A^a_jx^j$ to this system give rise to
a tableau $\A \subset\hom(V,W)$. Conversely, any tableau $\A
\subset\hom(V,W)$ determines a PDE system of this type.

\begin{remark}
$\A^{(q)}$ is the set of homogeneous polynomial solutions of degree
$q+1$ to \eqref{lin-pde}.
\end{remark}

\begin{theorem}
The {PDE} system associated to $\A$ is involutive $\iff$ $\A$ is
involutive.
\end{theorem}

The \textit{symbol} of $\eqref{lin-pde}$ is the annihilator $\B :=
\A^\perp \subset V\otimes W^\ast$ of $\A$. }
\end{example}

\vskip0.1cm
\begin{example}\label{pfaff-sys}
{\rm
Let $(\mathcal{I},\omega)$ be a Pfaffian differential system
(PDS) on a manifold $M$ with independence condition $\omega \neq 0$,
where
\[
  \mathcal{I} = \{\theta^1,\dots,\theta^s,d\theta^1,\dots,d\theta^s\}
 \quad \text{(algebraic ideal)}
   \]
and
$\omega = \omega^1\wedge \cdots\wedge \omega^n$.
Let $\pi^1,\dots,\pi^t$ be 1-forms such that
\[
 \theta^1,\dots,\theta^s;\, \omega^1,
  \dots,\omega^n;\, \pi^1,\dots,\pi^t
   \]
be a local adapted co\-fra\-me of $M$.

\vskip0.1cm

The Pfaffian differential system $(\mathcal{I},\omega)$ is called
\textbf{linear}\footnote{In the literature, other names are also
used to indicate linear systems: quasi-linear systems, systems in
good form, or systems in normal form.} if and only if
\[
 d\theta^a \equiv 0
\mod \{\theta^1,\dots,\theta^s, \omega^1,\dots,\omega^n\} \quad
(0\leq a \leq s).
  \]
The meaning of this condition is that the variety
$V_n(\mathcal{I},\omega)$ $\subset$ $G_n(TM,\omega)$ of integral
elements of $(\mathcal{I},\omega)$ is described by a system of
inhomogeneous linear equations (cf. Appendix \ref{ss:torsion}). A
linear PDS is described locally by
\[
 \left\{\begin{array}{l}
 \theta^a =0\\
 d\theta^a \equiv A^a_{\epsilon i}\pi^\epsilon \wedge\omega^i + \frac{1}{2}c^a_{i j}
 \omega^i \wedge\omega^j \mod \{\theta^1,\dots,\theta^s\}\\
 \omega = \omega^1\wedge \cdots\wedge \omega^n\neq 0,\\
 \end{array}\right.
    \]
where $c^a_{i j} = -c^a_{j i}$;\, $1\leq a \leq s$;\,
$1 \leq i,j\leq n$;\, $1\leq \epsilon \leq t$.

\vskip0.1cm Once we fix independent variables and take a point of
$M$, we can associate a tableau to the Pfaffian system as follows.
At $x\in M$, let $V^\ast = \span \{\omega^i\}$ and
$\{\frac{\partial}{\partial{\omega^i}}\}$ be the basis of its dual
$V$. Further, let $W^\ast = \span \{\theta^a \}$ and
$\{\frac{\partial}{\partial{\theta^a}}\}$ be the basis of its dual
$W$. We define a tableau $\A \subset W \otimes V^\ast$ by
\[
 \A := \span \{ A^a_{\epsilon i}\,\frac{\partial}{\partial{\theta^a}}\otimes
  {\omega^i}\, :\, \epsilon = 1,\dots,t\}.
  \]

The involutiveness of $(\mathcal{I},\omega)$ at $x$ is equivalent to
the involutiveness of $\A$ (algebraic condition) together with the
integrability condition $V(\mathcal{I},\omega)_{|x}$ $\neq
\emptyset$, which in turn is equivalent to the condition $c^a_{i j}
(x) = 0$, for each $a,i,j$ (``{torsion} vanishes at $x$''). See
Appendix \ref{ss:torsion} for more on the notion of torsion.

\begin{theorem}
The linear PDS $(\mathcal{I},\omega)$ is involutive at $x$
$\iff$ 1) $\A$ is involutive and 2) $V_n(\mathcal{I},\omega)_{|x} \neq \emptyset$.
\end{theorem}

The \textit{symbol} of $(\mathcal{I},\omega)$ is
the annihilator $\B := \A^\perp \subset V\otimes W^\ast$ of $\A$:
\[
  B = \span\{ B^\lambda = B^{\lambda i}_{a} \frac{\partial}{\partial{\omega^i}}\otimes
  {\theta^a} : B^{\lambda i}_{a} A^a_{\epsilon i}=0,\forall \lambda,\epsilon\}.
   \]
}
\end{example}

\section{Tableaux over Lie algebras}\label{s:lie-tableaux}

Let $(\g, [\,,])$ be a finite dimensional Lie algebra, $\a, \b$
vector subspaces of $\g$ such that $\a\oplus\b= \g$, and $\A \subset
\hom(\a,\b)$ a tableau. Define the polynomial map $\tau : \A \to
\b\otimes\Lambda^2(\a^\ast)$ by
\begin{multline}\nonumber
  \tau(Q)(A_1,A_2) := \left[A_1 + Q(A_1),A_2  + Q(A_2)\right]_{\b} \\
    - Q \left(\left[A_1 + Q(A_1),A_2 + Q(A_2)\right]_{\a}\right),
      \end{multline}
where $X_{\a}$ (resp. $X_{\b}$) denotes the $\a$ (resp. $\b$)
component of $X$.

\vskip0.1cm
\begin{definition}[\cite{MNtableaux}]\label{def-tableau}
A \textbf{tableau over} $\,\g$ is a tableau $\A \subset \hom(\a,\b)$
such that:
\begin{enumerate}
\item\label{1} $\A$ is involutive;
\vskip0.1cm
\item\label{2} $\tau(Q) \in \text{Im}\, \delta^{1,1} \subset \b \otimes
\Lambda^2 (\a^\ast)$, for each $Q\in \A$.
\end{enumerate}
\end{definition}

\vskip0.1cm
\begin{remark}
A detailed analysis of the examples at our disposal and
considerations about the problem of prolongation (cf. Remark
\ref{r:2-acyclic}) suggest that \textit{condition} (1) in the above
definition \textit{should be replaced by the condition that $\A$ is
2-acyclic}. As for condition (2) in the definition, it amounts to
the vanishing of a cohomology class in the group
\[
 H^{0,2}(\A) =
  \frac{\b \otimes \Lambda^{2} (\a^\ast)}{\delta^{1,1}
   \left(\A \otimes \a^\ast\right)}.
    \]
(cf. Remark \ref{r:2-acyclic} and Appendix \ref{ss:torsion} for the
the notion of torsion of a linear PDS).
\end{remark}

\vskip0.1cm
\begin{example}
{\rm If $\A \subset \hom (\a,\b)$ is involutive (or 2-acyclic) and
$\g$ is the abelian Lie algebra $\g =\a\oplus\b$, then $\tau(Q) =
0$, for each $Q\in \A$, and $\A$ can be considered as a tableau over
$\g$. Therefore, the concept of tableau over a Lie algebra is a
natural (non-commutative) generalization of the classical notion of
involutive (or 2-acyclic) tableau.}
\end{example}

\vskip0.1cm
\begin{example}
{\rm Let $\g$ be a semisimple Lie algebra with Killing form
$\langle\,, \rangle$. Let $\g=\g_0\oplus \g_1$ be a Cartan
decomposition. Then
\[
 [\g_0,\g_0]\subset \g_0, \quad [\g_0,\g_1]\subset \g_1, \quad [\g_1,\g_1]\subset \g_0.
  \]
Assume that $\mathrm{rank\,}\g/\g_0=k$ and that
$\a$ be a maximal ($k$-dimensional) abelian subspace of $\g_1$.
Then $\g_1 = \a\oplus \m$, where
\[
 \m = \a^\perp\cap \g_1
  \]
Further, let
\[
 \begin{array}{c}
 (\g_0)_\a = \left\{X\in \g_0 \,: \, [X,\a]=0\right\},\\
 (\g_0)_\a^\perp = \left\{X\in \g \,: \,
  \langle X, Y \rangle =0, \,\forall\, \, Y \in (\g_0)_\a \right\},\\
   \b :=\g_0\cap (\g_0)_\a^\perp.
    \end{array}
     \]

Then, for any regular element $A\in \a$, the maps
\[
 \text{ad}_A : \m \to \b, \quad \text{ad}_A : \b \to \m
  \]
are vector space isomorphisms and
\[
 X\in  \m \mapsto -\text{ad}_X \in \hom(\a,\b)
  \]
is injective, hence $\m$ \textit{can be identified with a linear subspace of}
$\hom(\a,\b)$.

\vskip0.1cm
\begin{proposition}[\cite{MNtableaux}]
If $\g$ is a semisimple Lie algebra and
$\a$, $\b$, and $\m$ are defined as above, then $\m$, regarded as a
subspace of $\hom(\a,\b)$, is a tableau over $\g$.
\end{proposition}

\vskip0.1cm
\begin{definition}
The tableau $\m$ is called a \textbf{Cartan tableau} over $\g$.
\end{definition}

\vskip0.1cm
\begin{remark}
As already indicated in the introduction, the idea of a tableau over
a Lie algebra has its origin in the method of moving frames and is
related to the existence of canonical adapted frames along generic
submanifolds in homogeneous spaces. The tableaux corresponding to
systems of submanifold geometry are constructed on the Lie algebras
of the transitive groups of transformations of the ambient spaces,
e.g., the Wilczynski--Cartan frame in projective differential
geometry (cf. \cite{MNtableaux}), or the canonical M\"obius frame in
conformal theory of surfaces (cf. Section \ref{ss:conf-ex}).
\end{remark}
}
\end{example}

\section{Differential systems associated with tableaux
over Lie algebras}\label{tableau-system}

Let $\A\subset \hom(\a,\b)$ be a tableau over $\g$ and let $G$ be
a connected Lie group with Lie algebra $\g$. We set $Y := G\times \A$
and refer to it as the \textit{configuration space}.

\vskip0.1cm
\begin{definition}
A basis $(A_1,\dots,A_k,B_1,\dots,B_h, C_1,\dots,C_s)$ of $\g$ is said \textit{adapted}
to $\A$ if
\begin{enumerate}
\vskip0.1cm
 \item $\a= \text{span}\,\{A_1,\dots,A_k\}$,
\vskip0.1cm
 \item $\mathrm{Im\,} \A := \sum_{Q\in \A} \mathrm{Im\,} Q =
\text{span}\,\{B_1,\dots,B_h\}$,
\vskip0.1cm
 \item $\b = \text{span}\,\{B_1,\dots,B_h,C_1,\dots,C_s\}$.
\end{enumerate}
\vskip0.2cm An adapted basis is {\it generic} if the flag
\[
 (0) \subset \text{span}\,\{ A_1\} \subset\cdots \subset \text{span}\,\{A_1,\dots,A_k\} =\a
  \]
is generic with respect to $\A$.

\end{definition}

For a generic adapted basis, let
\[
 (\alpha^1,\dots,\alpha^k,\beta^1,\dots,\beta^h,\gamma^1.\dots,\gamma^s)
  \]
denote the dual coframe on $G$. Given a basis
\[
 Q_\epsilon = Q^j_{\epsilon i} B_j \otimes \alpha^i \quad (\epsilon=1, \dots m)
  \]
of the tableau $\A$, $Y$ identifies with $G\times \R^m$ by
\[
 (g,  p^\epsilon Q_\epsilon) \in Y \mapsto (g;p^1,\dots,p^m) \in G\times \R^m.
  \]

\vskip0.1cm
\begin{definition} [\cite{MNtableaux}]
The EDS \textbf{associated with} $\A$ is the Pfaffian system
$(\mathcal{I}, \omega)$ on $Y$ generated (as a differential ideal)
by the linearly independent 1-forms
\[
 \left\{\begin{array}{l}
  \eta^j := \beta^j - p^\epsilon Q^j_{\epsilon i}\alpha^i, \quad (j=1,\dots,h),\\
   \gamma^1, \dots, \gamma^s, \\
    \end{array}\right.
     \]
with independent condition
$\omega = \alpha^1 \wedge \cdots\wedge \alpha^k \neq 0$.

\end{definition}

An immersed submanifold
\[
  \Phi = (g; p^1,\dots,p^m) : N^k \to G\times \A \cong G\times \R^m.
   \]
is an \textbf{integral manifold} of $(\I,\omega)$ if and only if
\begin{enumerate}
\vskip0.1cm
\item ${(\alpha^1 \wedge \cdots\wedge \alpha^k)}_{|N} \neq 0$;
\vskip0.1cm
\item $\beta^j = p^\epsilon Q^j_{\epsilon i}\alpha^i$, $j=1,\dots,h$;
\vskip0.1cm
\item $\gamma^1 = \cdots = \gamma ^s =0$.
\end{enumerate}

The main result in the construction is the following.

\vskip0.1cm
\begin{theorem} [\cite{MNtableaux}]\label{pds-inv}
Let $\A$ be a tableau over a Lie algebra $\g$. Then, the EDS
$(\I,\omega)$ associated with $\A$ is a linear PDS in involution. In
particular, the characters of $\A$ coincide with the Cartan
characters of $(\I,\omega)$.
\end{theorem}

\vskip0.1cm
\begin{remark}\label{r:2-acyclic}
Condition (2) in the definition of a tableau over a Lie algebra (cf.
Definition \ref{def-tableau}) tells us that the PDS associated with
a tableau over a Lie algebra is linear and with \textbf{vanishing
torsion} (cf. Appendix \ref{ss:torsion}). This together with the
condition that the tableau is 2-acyclic guarantee the existence of a
prolongation tower for the associated PDS which can be constructed
algebraically from the tableau and its prolongations. The
construction of the prolonged systems is a direct consequence of the
property of the tableau being 2-acyclic and is entirely based on the
Spencer cohomology of the tableau. The vanishing of the torsion is
needed only at the first step of the construction (cf. Remark
\ref{a:prolongation}). Therefore, the result stated in Theorem
\ref{pds-inv} can be generalized to the following.

\vskip0.1cm
\begin{theorem} [\cite{MNprep}]
Let $\A$ be a 2-acyclic tableau over a Lie algebra $\g$.
Then, the PDS $(\I,\omega)$ on $Y$ associated with $\A$ admits
regular prolongations of any order. Moreover, the construction of
prolongations is purely
algebraic. The configuration space of the $h$-prolonged system $(\I^{(h)},\omega)$
is
\[
  Y^{(h)} := G \times (\A\oplus\A^{(1)} \oplus \cdots \oplus \A^{(h)}).
     \]
If $k$ is the least integer such that $\A^{(k)}$ is involutive, then the
$k$-prolongation $(\I^{(k)},\omega)$ is in involution and its Cartan characters
coincide with that of $\A^{(k)}$.
\end{theorem}

\end{remark}

\section{Examples}\label{s:examples}

In this section, we illustrate the construction developed in Section
\ref{tableau-system} by discussing some examples.

\subsection{The PDS associated with an abelian tableau}
Let $\A \subset \hom(\R^k,\R^h)$ be an m-dimensional
involutive tableau over the (abelian) Lie algebra $\g = \R^k \oplus \R^h$,
spanned by the linearly independent ${h\times k}$ matrices
$Q_\epsilon = (Q^j_{\epsilon i})$.

We call $\A^{(1)}$-\textbf{system} the linear, homogeneous, constant
coefficient PDE system for maps
$P = (p^1,\dots,p^m) : \R^k \to \A\cong \R^m$ defined by the differential
inclusion $dP_{|x} \in \A^{(1)}$, for all $x\in \R^k$, where
$\A^{(1)}$ is the first prolongation of $\A$.
This system can be written  $$\delta^{1,1}(dP) = 0,$$ where $\delta^{1,1}$
is the Spencer coboundary of the tableau $\A$ (recall that
$\delta^{1,1} : C^{1,1} =\A \otimes (\R^k)^\ast \to C^{0,2}$).

\vskip0.1cm
\begin{lemma}
A map $P: \R^k \to \A\cong \R^m$ is a solution to the
$\A^{(1)}$-system if and only if the $\R^h$-valued 1-form
\[
 \theta = (\theta^1,\dots,\theta^h) \in \Omega^1(\R^k) \otimes \R^h,
 \quad \theta^j = p^\epsilon Q_{\epsilon a}^j dx^a,
    \]
is closed.
\end{lemma}

As a consequence, we have

\vskip0.1cm
\begin{corollary}
Let $P: \R^k \to \A\cong \R^m$ be a solution to the
$\A^{(1)}$-system and let $y = (y^1,\dots,y^h)$ be a primitive of $\theta$
(i.e., $\theta = dy$).
Then
$$
 \R^k \ni x \mapsto (x,y(x),P(x)) \in \R^k\oplus \R^h \oplus \R^m
  $$
is an integral manifold of the PDS $(\I,\omega)$ associated with $\A$.
Moreover, every integral manifold of $(\I,\omega)$ arises in this way.
\end{corollary}

We can conclude that the integral manifolds of $(\I,\omega)$ correspond
to the solutions of the $\A^{(1)}$-system. Moreover, $(\I,\omega)$ is
in involution (as a differential system) and the Cartan characters coincide
with those of the tableau $\A$.

\subsection{The PDS associated with a Cartan tableau and the $G/G_0$-system}
Let $G/G_0$ be a semisimple symmetric space of rank $k$ and $\g = \g_0 \oplus \g_1$
a Cartan decomposition of $\g$.
Let $(A_1,\dots,A_k)$ be a regular basis for the maximal
abelian subalgebra $\a \subset \g_1$.
According to Terng \cite{Ter97}, the $G/G_0$-\textbf{system}
(or the $k$-dimensional system associated to $G/G_0$)
is the system of PDEs for maps $F : U\subset \a \to \m :=\g_1\cap \a^\perp\,$ defined by
\[
 \left[A_i, \frac{\partial F}{\partial x^j}\right]-
  \left[A_j, \frac{\partial F}{\partial x^i}\right]= \left[[A_i, F], [A_j, F]\right],
  \]
$1\leq i<j\leq k$, where $x^i$ are the coordinates
with respect to $(A_1,\dots,A_k)$.

\vskip0.1cm
\begin{lemma}
A map $F : \a \to \m$ is a solution of the $G/G_0$-system
if and only if the $\g$-valued 1-form
\[
  \theta = \alpha + [\alpha,F] \in \Omega^1(\a)\otimes \g,
  \]
satisfies the Maurer--Cartan equation $d\theta +
\frac{1}{2}[\theta\wedge\theta] =0$, where $\alpha = \alpha^i\otimes
A_i$ is the tautological 1-form on $\a$.
\end{lemma}

\vskip0.1cm
\begin{corollary}
 Let $F : \a \to \m$ be a solution of the
$G/G_0$-system and let $g : \a \to G$ be a \textit{primitive} of $\theta$
(i.e. a solution to $g^{-1}dg = \theta$). Then
$$
   \a \ni x \mapsto (g(x),F(x)) \in G\times \m
  $$
is an integral manifold of the PDS $(\mathcal{I}, \omega)$ on $Y = G \times
\m$ associated with the Cartan tableau $\m \subset \hom(\a,\b)$.
Conversely, any integral manifold of $(\mathcal{I}, \omega)$ arises in this way.
\end{corollary}

In conclusion, the integral manifolds of the PDS $(\I,\omega)$
associated with the Cartan tableau $\m \subset \hom(\a,\b)$ are given
by the solutions of the corresponding $G/G_0$-system. Moreover, $(\I,\omega)$ is
in involution and its Cartan characters coincide with those
of the tableau $\m$ (i.e., $s_1 = n$, $s_j =0$, $j = 2,\dots,n$).
In particular, the
general solution depends on $n$ functions in one variable.

\vskip0.1cm
\begin{remark}
If $(A_1,\dots,A_n)$ is a basis of $\a$,
$(x^1,\dots,x^n)$ the corresponding coordinates, and
$F=(F^1,\dots,F^{m-n}) : U \subset \a \to \b$, then the
$G/G_0$-system can be written
\[
 B^a_{\alpha,i}\frac{\partial F^{\alpha}}{\partial x^j}-B^a_{\alpha,j}
  \frac{\partial F^{\alpha}}{\partial x^i} =\Phi^a_{ij},
    \]
where the coefficients $B^a_{\alpha i}$ are constant and $\Phi_{ij}$
are analytic functions.

In general, the PDS associated to a tableau over a Lie algebra corresponds to
the nonlinear system of equations
\[
 B^a_{\alpha,i}\frac{\partial F^{\alpha}}{\partial x^j}-
  B^a_{\alpha,j}\frac{\partial F^{\alpha}}{\partial x^i}
  +B^a_{\alpha,\beta}\left(\frac{\partial F^{\alpha}}{\partial x^i}
   \frac{\partial F^{\beta}}{\partial x^j}-
   \frac{\partial F^{\alpha}}{\partial x^j}
    \frac{\partial F^{\beta}}{\partial x^i}\right)=\Phi^a_{ij}.
     \]
\end{remark}

\subsection{Old and new involutive systems in conformal surface
theory}\label{ss:conf-ex}

In this section we discuss some old and new involutive
systems/tableaux arising in conformal geometry of surfaces. We start
by recalling some preliminary material. Consider Minkowski 5-space
${\R}^{4,1}$ with linear coordinates $x^0,\dots,x^4$ and Lorentz
scalar product given by the quadratic form
\begin{equation}
 \langle x,x\rangle = -x^0x^4
   +(x^1)^2+(x^2)^2+(x^3)^2=\eta_{ij}x^{i}x^{j}.\label{product}
      \end{equation}
Classically, the M\"obius space ${S}^3$ (conformal 3-sphere) is
realized as the projective quadric $\left\{[x]\in {\RP}^4 : \langle
x,x\rangle =0 \right\}$. Accordingly, ${S}^3$ inherits a natural
conformal structure and the identity component $G \cong \SO_0(4,1)$
of the pseudo-orthogonal group of (\ref{product}) acts transitively
on ${S}^3$ as group of orientation-preserving, conformal
transformations (see \cite{BrWill}). The Maurer--Cartan form of $G$
will be denoted by $\omega = (\omega^i_j)$.

Let $f : {{U}}\subset{\R}^2 \to {{S}}^3$ be an umbilic free,
conformal immersion. A \textit{M\"obius frame field} along $f$ is a
map ${{g}}=(g_0,\dots,g_4) : {{U}} \to G$ such that $f(p)=[g_0(p)]$,
for all $p\in {U}$. According to \cite{BrWill}, there exists a
canonical M\"obius frame field\footnote{If ${g}$ is a canonical
frame, any other canonical frame is given by
$(g_0,-g_1$,$-g_2,g_3,g_4)$.} ${{g}} : {{U}} \to G$ along $f$ such
that its Maurer--Cartan form $\beta =(\beta^i_j)=
{{g}}^{\ast}\omega$ takes the form
$$
\left(\negthickspace\begin{array}{ccccc}
-2q_2\beta^1_0+2q_1\beta^2_0&p_1\beta^1_0 +p_2\beta^2_0
&-p_2\beta^1_0+
p_3\beta^2_0&0&0\\
\beta^1_0&0&-q_1\beta^1_0-q_2\beta^2_0&-\beta^1_0&p_1\beta^1_0 +p_2\beta^2_0 \\
\beta^2_0&q_1\beta^1_0+q_2\beta^2_0&0&\beta^2_0&-p_2\beta^1_0+ p_3\beta^2_0\\
0&\beta^1_0&-\beta^2_0&0&0\\
0&\beta^1_0&\beta^2_0&0&2q_2\beta^1_0-2q_1\beta^2_0
\end{array} \negthickspace\right)\label{centralframe}
$$
with $\beta^1_0\wedge \beta^2_0 > 0$. The smooth functions $q_1$,
$q_2$, $p_1$, $p_2$, $p_3$ form a complete system of conformal
invariants for $f$ and satisfy the following structure equations
\begin{eqnarray}
d\beta^1_0=-q_1\beta^1_0\wedge\beta^2_0,&{} &
d\beta^2_0=-q_2\beta^1_0\wedge\beta^2_0,\label{cf1} \\
dq_1\wedge\beta^1_0+dq_2\wedge\beta^2_0&=&(1+p_1+p_3+{q_1}^2+{q_2}^2)\beta^1_0\wedge
\beta^2_0,\label{cf2} \\
dq_2\wedge\beta^1_0-dq_1\wedge\beta^2_0&=& - p_2\beta^1_0\wedge\beta^2_0,\label{cf3}\\
dp_1\wedge\beta^1_0 + dp_2\wedge \beta^2_0&=&(4q_2p_2
+q_1(3p_1+p_3))\beta^1_0\wedge\beta^2_0, \label{cf4}\\
dp_2 \wedge\beta^1_0 - dp_3\wedge\beta^2_0&=&(4q_1p_2
-q_2(p_1+3p_3))\beta^1_0\wedge\beta^2_0.\label{cf5}
\end{eqnarray}

\vskip0.2cm

\subsubsection{The M\"obius tableau} The existence of a canonical
M\"obius frame field suggests the following construction. Let
$(\alpha^1$, $\alpha^2$, $\beta^1$, $\dots$, $\beta^4$, $\gamma^1$,
  $\dots$, $\gamma^4)$ be the basis of $\g^\ast$ defined by
$$
\begin{cases}
 \alpha^1 = \omega^1_0, \quad \alpha^2 = \omega^2_0, \quad \beta^1 =
  \omega^0_0, \quad \beta^2 = \omega^0_1, \quad \beta^3 = \omega^0_2,
   \quad \beta^4 = \omega^2_1,\\
   \gamma^1 = \omega^0_3, \quad \gamma^2 = \omega^3_0, \quad
    \gamma^3 =  \omega^1_0 -\omega^3_1, \quad \gamma^4 =  \omega^2_0 +
     \omega^3_2.
\end{cases}
$$
Next, let $(A_1, A_2, B_1,\dots,B_4, C_1,\dots,C_4)$ be its dual
basis and set
$$
 \a = \span \{A_1, A_2\}, \quad \b = \span \{B_1,\dots, B_4\}.
  $$
Consider the 5-dimensional subspace $\M\subset \hom(\a,\b)$
consisting of all elements $Q(q,p)$ of the form
\begin{multline}\nonumber
  Q(q,p) = q_1(B_4\otimes \alpha^1 + 2B_1 \otimes \alpha^2) +
            q_2 (-2B_1\otimes \alpha^1 + B_4 \otimes \alpha^2)\\
            + p_1B_2\otimes \alpha^1 + p_2(-B_3\otimes \alpha^1 + B_2 \otimes
            \alpha^2) + p_3B_3\otimes \alpha^2,
               \end{multline}
where $q=(q_1,q_2)\in \R^2$, $p=(p_1.p_2,p_3) \in \R^3$. A direct
computation yields the following.

\vskip0.2cm
\begin{lemma}
The subspace $\M$ is a tableau over $\g\cong \so(4,1)$.
\end{lemma}

\vskip0.2cm

The PDS associated with the tableau $\M$, referred to as the
M\"obius system, is the PDS on $Y=G\times \M\cong G\times \R^5$
generated by the 1-forms $\gamma^1,\dots, \gamma^4, \eta^1, \dots,
\eta^4$, where
$$
\begin{cases}
 \eta^1 = \beta^1 + 2q_2\omega^1_0 -2q_1\omega^2_0, \quad \eta^2 = \beta^2 -p_1\omega^1_0
 -p_2\omega^2_0,\\
 \eta^3 = \beta^3 + p_2\omega^1_0 -p_3\omega^2_0,
 \quad \eta^4 = \beta^4 -q_1 \omega^1_0 -q_2\omega^2_0,
   \end{cases}
     $$
with independence condition $\omega^1_0\wedge\omega^2_0\neq 0$. The
integral manifolds of the M\"obius system are the 2-dimensional
submanifolds
$$
 (g;q,p) : M^2 \to G \times \M \cong G \times \R^5
   $$
such that:
\begin{itemize}
 \item $f=[g_0] \to S^3$ is an umbilic free, conformal immersion;
 \item $g : M^2 \to G$ is a canonical M\"obius frame along $f$;
 \item $q_1$, $q_2$, $p_1$, $p_2$, $p_3 : M^2 \to \R$ are the conformal
    invariants of $f$.
     \end{itemize}

\subsubsection{Willmore surfaces}

\textit{Willmore immersions} are defined as extremals for the
Willmore functional $\int{(H^2-K)}dA$ ($H$ the mean curvature, $K$
the Gauss curvature). They are characterized by the Euler--Lagrange
equation
$$
 \Delta{H}+2H(H^2-K)=0,
   $$
which expressed in terms of the conformal invariants is equivalent
to the equation $p_1 = p_3$ (cf. \cite{BrWill, MNvarsavia,
Willmore}). Willmore surfaces can be seen as integral manifolds of
the M\"obius system restricted to the submanifold of $Y$ given by
$$
 Y_{W}= \{Q(q,p) \in Y \, | \, p_1 = p_3\}.
  $$
Now, it is easy to check that the subspace
$$
 \M_W  := \{Q(q,p) \in \M \, | \, p_1 = p_3\}
  $$
defines a 4-dimensional tableau over $\g\cong \so(4,1)$ with
characters $s_0=8$, $s_1= 4$, $s_2 =0$, and that $Y_{W}$ is the
configuration space of $\M_W$. Observe also that the restriction to
$Y_{W}$ of the M\"obius system is exactly the PDS associated with
$\M_W$.

\subsubsection{Other classes of surfaces}
More generally, one could consider the class of surfaces whose
invariant functions $p_1$ and $p_3$ satisfy a linear relation, that
is, are expressed by $p_1(t) = t\cos a +b_1$, $p_2(t) = t\sin a
+b_2$, for real constants $a$, $b_1$, $b_2$. This class includes
Willmore surfaces as special examples and corresponds to the
4-dimensional affine tableau
$$\M_{(a,b_1,b_2)} =\{Q(q,p_1(t),p_2,p_3(t)) \, | \, t, p_2\in\R,
q\in \R^2\}.$$ Also in this case, an algebraic, direct computation
shows that $\M_{(a,b_1,b_2)}$ is involutive. Therefore, by the
construction developed in the previous section, the associated PDS
is in involution. Its Cartan characters are $s_0=8$, $s_1= 4$, $s_2
=0$.

\vskip0.2cm
\begin{remark}
Similar arguments have been used in connection with the study of
surfaces in projective differential geometry \cite{MNtableaux}. The
same approach can also be used to discuss surface theory in the
framework of Laguerre geometry \cite{MNtransactions} and other
classical geometries.
\end{remark}

\section{Further developments}\label{s:developments}

\subsection{}
Continue the program, initiated with the study of several classes of
integrable surfaces in projective differential geometry
\cite{MNtableaux}, of identifying the geometry associated to a given
tableau/system, i.e., find submanifolds in some homogeneous space
whose integrability conditions are given by the PDS associated with
the given tableau.

\subsection{}
Study the algebraic structure of tableaux over Lie algebras to
understand when a tableau generates an integrable geometry.

\subsection{}
Study the Cauchy problem for the associated systems (cf.
\cite{MNjmp, MNphysD}).

\subsection{}
Analyze the characteristic cohomology of a tableau over a Lie
algebra, its geometric interpretations, and its relations with the
characteristic cohomology of Bryant--Griffiths \cite{BG1, BG2}.

\appendix

\section{The Spencer complex}\label{ss:spencer}[cf. \cite{BCGGG}]

Retaining the notation of Section \ref{s:tableaux}, identify the
symmetric product $S^{q}(\a^\ast)$ with the space of homogeneous
polynomials of degree $q$ on $\a$. For each $v\in \a$, let
$\delta_v$ be the map of $\b \otimes S^{q}(\a^\ast) \to \b \otimes
S^{q-1}(\a^\ast)$ given by partial differentiation w.r.t. $v$. Let
$v_1,\dots,v_n$ be a basis of $\a$, and $v^1,\dots,v^n$ its dual
basis.

The operator
\[
  \b\otimes S^{q}(\a^\ast) \otimes \Lambda^p(\a^\ast) \xrightarrow{\delta^{q,p}}
    \b\otimes S^{q-1}(\a^\ast) \otimes \Lambda^{p+1}(\a^\ast)
      \]
given by
\[
 \delta^{q,p}\xi : = \sum \delta_{v_i}\xi\wedge v^i
   \]
($\delta^{0,p} =0$, for $p\geq 0$) is independent of the basis,
$\delta^2 = 0$, and the sequence of the corresponding bigraded
complex is exact except when $q=0$ and $p=0$.

Let $\A\subset \hom(\a,\b)$ be a tableau with
prolongations ${\A}^{(h)}$, $h\geq 0$.
Consider the sequence of spaces
\[
 C^{q,p}(\A):={\A}^{(q-1)}\otimes \Lambda^p(\a^\ast),
 \]
for integers $q\geq 0$ and $0 \le p \leq n$.
Note that since ${\A}^{(q-1)} \subset \b \otimes S^{q}(\a^\ast)$, the space
$C^{q,p}(\mathbf{A})$ is a subspace of
$\b\otimes S^{q}(\a^\ast) \otimes \Lambda^p(\a^\ast)$.
We have
\[
  \delta C^{q,p}(\A) \subset C^{q-1,p+1}(\A),
    \]
but the sequence
\[
 C^{q+1,p-1}(\A) \xrightarrow{\delta^{q+1,p-1}} C^{q,p}(\A)
   \xrightarrow{\delta^{q,p}} C^{q-1,p+1}(\A)
   \]
is no longer exact for all $p$ and $q$.
The associated cohomology groups
\[
  H^{q,p}(\A) := Z^{q,p}(\A)/B^{q,p}(\A)
   \]
are called the \textit{Spencer groups} of $\A$,
where $B^{q,p}(\A)=\mathrm{Im}(\delta^{q+1,p-1})$ and
$Z^{q,p}(\A)=\ker(\delta^{q,p})$.
Notice that
$Z^{0,p}(\A)=\b\otimes\Lambda^p(\a^\ast)$, for all $p\geq 0$, and
   $Z^{q,1}(\A)={\A}^{(q)}$, for all $q\geq 1$.

\vskip0.1cm
A significant result in the subject is that the vanishing of the
$H^{q,p}$ is equivalent to involutiveness.

\begin{theorem}[\cite{GuiSte64}]
A tableau $\A$ is involutive if and only if $H^{q,p}(\A)$ is zero,
for all $q\geq 1$ and $p\geq 0$.
\end{theorem}

A weaker condition than involutiveness is the following.

\begin{definition}
A tableau $\A$ is called 2-\textit{acyclic} if $H^{q,2}(\A)=(0)$, for all $q\ge 1$.
\end{definition}

Another way of formulating the condition
\[
 H^{q,p}(\A) = (0), \quad \text{for all} \quad q\geq 1,\, p\geq 0
  \]
is that the sequences
\[
 0\rightarrow \A^{(k)}\xrightarrow{\delta} \A^{(k-1)}\otimes \a^\ast \rightarrow
 \cdots
  \xrightarrow{\delta} \A \otimes\Lambda^{k-1} (\a^\ast) \rightarrow
  \]
\[
 \cdots\xrightarrow{\delta} \b \otimes \Lambda^{k} (\a^\ast) \rightarrow
  \frac{\b \otimes \Lambda^{k} (\a^\ast)}{\delta \left(\A \otimes\Lambda^{k-1} (\a^\ast)\right)}
   \rightarrow 0
    \]
are exact for all $k$. In particular, we have
\[
 H^{0,k}(\A) = \frac{\ker(\delta^{0,k})}{\im(\delta^{1,k-1})}
  =\frac{\b \otimes \Lambda^{k} (\a^\ast)}{\delta
   \left(\A \otimes\Lambda^{k-1} (\a^\ast)\right)}.
    \]

\section{Torsion of a Pfaffian systems}\label{ss:torsion}[cf. \cite{BCGGG}]

Retaining the notation of Example \ref{pfaff-sys},
let $(\mathcal{I},\omega)$ be a Pfaffian system.
An admissible integral element $E\in V_n(\mathcal{I},\omega)_{|x}$ is given by
\[
 \theta^a = 0, \quad \pi^\epsilon = p^\epsilon_i \omega ^i,
  \]
where the fiber coordinates $p^\epsilon_i$ satisfy
\[
 A^a_{\epsilon \, j}(x)p^\epsilon_i - A^a_{\epsilon \, i}(x)p^\epsilon_j
   + c^a_{i\, j}(x) = 0.
    \]
Under a change of coframe
\begin{equation}\label{coframe-change}
 \tilde{\theta}^a = \theta^a, \quad \tilde{\omega}^i=\omega^i, \quad
\tilde{\pi}^\epsilon = \pi^\epsilon - p^\epsilon_i\omega^i
  \end{equation}
the numbers $c^a_{i\, j}(x)$ transform to
\[
 \tilde{c}^a_{ij}(x) = A^a_{\epsilon \, j}(x)p^\epsilon_i -
  A^a_{\epsilon \, i}(x)p^\epsilon_j + c^a_{i\, j}(x).
   \]
This defines an equivalence relation
\[
 \tilde{c}^a_{ij}(x)\sim {c}^a_{ij}(x).
  \]

\begin{definition}
 The equivalence class $[{c}^a_{ij}(x)]$ is called the
\textbf{torsion} of $(\mathcal{I},\omega)$.
\end{definition}

\vskip0.1cm
\begin{lemma}
 The torsion of $(\mathcal{I},\omega)$ lives in
\[
 H^{0,2}(\A) = \frac{W \otimes \Lambda^{2} (V^\ast)}{\delta^{1,1}
  \left(\A \otimes V^\ast\right)} = \frac{\ker (\delta^{0,2})}{\im (\delta^{1,1})}.
    \]
\end{lemma}

\begin{proof} If
$Q =
  p^\epsilon_j A^a_{\epsilon \, i}\,\frac{\partial}{\partial{\theta^a}}\otimes
  {\omega^i} \otimes {\omega^j} \in \A \otimes V^\ast$, then
\[
 \delta^{1,1} \left(Q \right) = \sum_{i< j}\left(p^\epsilon_j A^a_{\epsilon \, i} -
   p^\epsilon_i A^a_{\epsilon \, j}\right) \,\frac{\partial}{\partial{\theta^a}}\otimes
    {\omega^i} \wedge {\omega^j}.
     \]
According to the transformation rule \eqref{coframe-change} of the ${c}^a_{ij}$ under a
coframe change, the cocycle
\[
 \frac{1}{2}c^a_{i\, j}(x)\,\frac{\partial}{\partial{\theta^a}}\otimes
  {\omega^i} \wedge {\omega^j} \in C^{0,2}(\A)
   \]
gives a class in $H^{0,2}(\A)$.
\end{proof}

\vskip0.1cm
\begin{remark}
The vanishing of the torsion is a necessary
and sufficient condition for the existence of an integral element over $x$.
\end{remark}

\vskip0.1cm
\begin{lemma}
The torsion of $(\mathcal{I}^{(1)}, \omega)$ lives in the
vector spaces
\[
  H^{0,2}(\A^{(1)}) \cong H^{1,2}(\A).
   \]
\end{lemma}

We also recall the following.

\begin{lemma}
\[
 H^{q,p}(\A^{(1)}) \cong H^{q+1,p}(\A), \quad q\geq 1.
  \]
   \end{lemma}

\begin{remark}\label{a:prolongation}
Thus, the involutiveness of the tableau $A$ associated to
$(\mathcal{I}, \omega)$ implies both the involutiveness of the
tableau $A^{(1)}$ and the vanishing of torsion of the prolonged
system $(\mathcal{I}^{(1)}, \omega)$.
\end{remark}

\end{document}